\documentclass[11pt]{article}
\usepackage[a4paper,margin=1in]{geometry}
\usepackage{amsmath,amssymb,amsthm,mathtools}
\usepackage{microtype}
\usepackage[T1]{fontenc}

\title{Dense signed sums of non-integer powers}
\author{David Treeby}
\date{}

\theoremstyle{plain}
\newtheorem{theorem}{Theorem}
\newtheorem{proposition}[theorem]{Proposition}
\newtheorem{lemma}[theorem]{Lemma}

\theoremstyle{definition}

\newtheorem{remark}[theorem]{Remark}

\begin{document}
\maketitle

\section{Introduction}

A familiar feature of conditional convergence is that a sequence of small
terms may be signed so as to approach essentially any prescribed value. More
precisely, if \(a_k>0\), \(a_k\to0\), and \(\sum_{k=1}^{\infty}a_k=\infty\),
then for every \(r\in\mathbb R\) there are signs
\(\varepsilon_k\in\{\pm1\}\) such that
\[
\sum_{k=1}^{N}\varepsilon_k a_k \longrightarrow r .
\]
Thus small but non-summable increments can be steered towards any target. 

The purpose of this note is to show that an analogous phenomenon can occur \emph{even when the underlying terms $a_k$ do not tend to zero}. Indeed, it can occur for sequences whose terms grow without bound. For example, take the sequence
\[1^{5/2},2^{5/2},3^{5/2},4^{5/2},5^{5/2},\ldots\]
Though these terms grow rapidly, we shall show that there is a choice of signs \(\varepsilon_k\in\{\pm1\}\) for which the set of partial sums
\(
S_N=\sum_{k=1}^{N}\varepsilon_k k^{5/2}
\)
is dense in \(\mathbb R\). Our main aim is to prove the following. 

\begin{theorem}[Main Theorem]\label{thm:main-theorem}
Let \(j>0\) be non-integer. Then there exists a sequence of signs
\(\varepsilon_k\in\{\pm1\}\) such that the partial sums
\[
S_N=\sum_{k=1}^{N}\varepsilon_k k^j
\]
are dense in \(\mathbb R\). Equivalently, for every \(r\in\mathbb R\), there is
an increasing sequence \(N_1<N_2<\cdots\) such that
\(
S_{N_\ell}\longrightarrow r .
\)
\end{theorem}

\noindent  This extends the work of Treeby and Wang, who showed that signed sums of non-integer powers of consecutive integers can be made arbitrarily close to zero \cite{TreebyWang}. Although this theorem seems implausible, the key observation is that we can create high-order cancellation among consecutive values. In the non-integer setting these cancellations are not exact, but they produce a sequence of small block sums whose magnitudes tend to $0$ while their total variation is still infinite. Once such a sequence of ``available increments'' is in hand, a classical greedy argument steers the partial sums to any prescribed limit. 

For integer powers, there is a corresponding discrete result. Bleicher \cite{Bleicher} proved that for every positive integer \(m\) and every \(z\in\mathbb Z\), there exist
\(N\ge1\) and signs \(\varepsilon_k\in\{\pm1\}\) such that
\[
\sum_{k=1}^{N}\varepsilon_k k^m=z.
\]
Yu \cite{Yu} later extended this to integer-valued polynomials with no fixed divisor greater than \(1\). In particular, one may choose a single sequence of signs for which the set of partial sums is exactly \(\mathbb Z\). Thus integer powers fill \(\mathbb Z\), while the non-integer powers considered here densely fill \(\mathbb R\).

\section{A steering lemma}\label{sec:steer}

The following is the basic ``crossing'' principle. It is closely related to the
usual rearrangement phenomenon for conditionally convergent series, but is
phrased in the signing language used throughout this note.

\begin{lemma}\label{lem:steer}
Let \(a_k>0\), \(a_k\to0\), and \(\sum_{k=1}^{\infty}a_k=\infty\). Then for
every \(r\in\mathbb R\), there are signs
\(\varepsilon_k\in\{\pm1\}\) such that
\[
S_N=\sum_{k=1}^{N}\varepsilon_k a_k \longrightarrow r .
\]
\end{lemma}

\begin{proof}
Set \(S_0=0\). Having chosen \(\varepsilon_1,\dots,\varepsilon_{k-1}\), choose
\[
\varepsilon_k=
\begin{cases}
+1,& S_{k-1}\le r,\\
-1,& S_{k-1}>r.
\end{cases}
\]
Thus each step moves the partial sum towards \(r\).

We claim that the sequence crosses \(r\) infinitely often. Indeed, if it
eventually stayed below \(r\), then every subsequent sign would be \(+1\), and
the divergence of \(\sum a_k\) would eventually force the partial sums above
\(r\). The case where the sequence eventually stayed above \(r\) is identical.
Thus crossings occur arbitrarily far along the sequence.

Let \(\epsilon>0\). Since \(a_k\to0\), there is \(N_0\) such that \(a_k<\epsilon\) for
all \(k\ge N_0\). Choose a crossing time \(N\ge N_0\). Then
\[
|S_N-r|<a_N<\epsilon.
\]
Moreover, once \(|S_N-r|<\epsilon\) and all later steps have size less than
\(\epsilon\), the greedy rule keeps every subsequent partial sum within \(\epsilon\) of
\(r\). Hence \(S_N\to r\).
\end{proof}

\section{Thue--Morse blocks and power asymptotics}\label{sec:block}

Fix an integer \(n\ge1\). Let
\(s_1,\dots,s_{2^n}\in\{\pm1\}\) be the length-\(2^n\) Thue--Morse sign pattern,
defined by
\[
s_k=(-1)^{\nu(k-1)},
\]
where \(\nu(m)\) denotes the number of \(1\)'s in the binary expansion of
\(m\). Thus, for example, when \(n=3\), the sign pattern is
\[
+,\,-,\,-,\,+,\,-,\,+,\,+,\,-.
\]
\noindent The associated generating polynomial is
\[
\sum_{k=1}^{2^n}s_k z^{k-1}
=
\prod_{q=0}^{n-1}\left(1-z^{2^q}\right).
\]
Since the right-hand side has a zero of order \(n\) at \(z=1\), the first
\(n\) moments vanish. That is, for \(0\leq i\leq n-1\),
\begin{equation}\label{eq:moments}
\sum_{k=1}^{2^n}s_k k^i=0. 
\end{equation}
The next moment does not vanish; in particular,
\begin{equation}\label{eq:next-moment}
\sum_{k=1}^{2^n}s_k k^n\ne0.
\end{equation}
This is the usual Prouhet--Tarry--Escott cancellation encoded by the Thue--Morse signs \cite{AlloucheShallit}. Moreover, it gives
\(\sum_{k=1}^{2^n}s_k P(k)=0\) for every polynomial \(P\) of degree less than
\(n\).

We now use these signs to form signed blocks of consecutive values of a function. The subscript \(n\) records that we are using the length-\(2^n\) Thue--Morse block. If \(f\) is a real-valued function defined for sufficiently large positive \(x\), define
\begin{equation}\label{eq:blockop}
(B_n f)(x)=\sum_{k=1}^{2^n}s_k f(x+k).
\end{equation}
Thus \(B_n f(x)\) is the signed contribution of the block
\[
f(x+1),f(x+2),\dots,f(x+2^n)
\]
using the fixed length-\(2^n\) Thue--Morse signs \(s_1,\dots,s_{2^n}\). In particular, if
\(a_k=f(k)\), then \(B_n f(m)\) is the signed block sum formed from
\[
a_{m+1},a_{m+2},\dots,a_{m+2^n}.
\]
For the power function \(f(x)=x^j\), this becomes
\[
(B_n f)(x)=\sum_{k=1}^{2^n}s_k(x+k)^j.
\]

\subsection*{Exact cancellation for integer powers}

If \(f(x)=x^m\), where \(m\) is a nonnegative integer with \(m<n\), then
\((B_n f)(x)\equiv0\). Indeed, for fixed \(x\), the expression \((x+k)^m\)
is a polynomial in \(k\) of degree less than \(n\), so the cancellation
property above gives
\[
(B_n f)(x)=\sum_{k=1}^{2^n}s_k(x+k)^m=0.
\]
This is the classical Prouhet--Tarry--Escott phenomenon encoded in
\eqref{eq:moments}.

\subsection*{Leading asymptotics for non-integer powers}

The following estimate is the basic non-integer analogue of the exact cancellation above. It is essentially the block asymptotic underlying the work of Treeby and Wang on non-integer powers \cite{TreebyWang}; we include the short proof here as it exhibits the cancellation mechanism in action. 

Let \(j>0\) be non-integer and choose
\[
n=\lceil j\rceil,
\qquad\text{so that}\qquad n-1<j<n.
\]
For the power function \(f_j(x)=x^j\), write
\[
F_j(x)=(B_n f_j)(x)=\sum_{k=1}^{2^n}s_k(x+k)^j.
\]
The moment vanishing \eqref{eq:moments} implies that the first non-vanishing
term in the binomial expansion occurs at order \(n\).

\begin{proposition}[Power block asymptotic]\label{prop:asymptotic}
Let \(j>0\) be non-integer and \(n=\lceil j\rceil\). Then there is a nonzero
constant \(C_{j,n}\) such that
\[
F_j(x)=C_{j,n}x^{j-n}+O(x^{j-n-1})
\qquad (x\to\infty).
\]
In particular, since \(j-n\in(-1,0)\), the block values \(F_j(x)\) tend to \(0\), but their absolute values are not summable along any arithmetic progression.
\end{proposition}

\begin{proof}
For large \(x\), expand
\[
(x+k)^j=x^j\left(1+\frac{k}{x}\right)^j
=x^j\sum_{i=0}^{n+1}\binom{j}{i}\frac{k^i}{x^i}
+O(x^{j-n-2}),
\]
where the error term is uniform for \(1\le k\le2^n\). Summing with the
Thue--Morse signs gives
\[
F_j(x)
=x^j\sum_{i=0}^{n+1}\binom{j}{i}\frac{1}{x^i}
\left(\sum_{k=1}^{2^n}s_k k^i\right)
+O(x^{j-n-2}).
\]
The inner moment vanishes for \(i<n\) by \eqref{eq:moments}. Hence the leading
term is the \(i=n\) term:
\[
F_j(x)
=
x^{j-n}\binom{j}{n}
\left(\sum_{k=1}^{2^n}s_k k^n\right)
+O(x^{j-n-1}).
\]
Thus
\[
C_{j,n}=\binom{j}{n}\sum_{k=1}^{2^n}s_k k^n.
\]
This constant is nonzero because \(\binom{j}{n}\ne0\) for non-integer \(j\),
and the \(n\)-th moment is nonzero by \eqref{eq:next-moment}.
\end{proof}

\begin{remark}
For the present argument, only the non-vanishing of \(C_{j,n}\) matters. The
precise value of the moment
\(\sum_{k=1}^{2^n}s_k k^n\) can also be computed explicitly from the
finite-difference interpretation of the Thue--Morse block.
\end{remark}

\section{Proof of the main theorem}\label{sec:main}

We now convert the block asymptotic into a single choice of signs whose partial sums are dense in \(\mathbb R\).

\begin{proof}[Proof of Theorem \ref{thm:main-theorem}]
Let \(j>0\) be non-integer and set \(n=\lceil j\rceil\). For each block index \(t\ge0\), define
\[
b_t=F_j(t2^n)=\sum_{k=1}^{2^n}s_k(t2^n+k)^j.
\]
By Proposition \ref{prop:asymptotic},
\[
b_t=C_{j,n}(t2^n)^{j-n}+O(t^{j-n-1})
\qquad (t\to\infty),
\]
where \(C_{j,n}\ne0\). Since \(j-n\in(-1,0)\), it follows that, after discarding at most finitely many initial blocks,
\[
|b_t|\to0
\qquad\text{and}\qquad
\sum |b_t|=\infty.
\]

Choose the signs in the discarded initial blocks arbitrarily. It remains to choose, for each sufficiently large block, either the Thue--Morse pattern
\[
s_1,\dots,s_{2^n}
\]
or its negative. That is, for each remaining block choose a sign \(\sigma_t\in\{\pm1\}\), and set
\[
\varepsilon_{t2^n+k}=\sigma_t s_k
\qquad (k=1,\dots,2^n).
\]
The contribution of this block to the partial sum is then
\[
\sum_{k=1}^{2^n}\varepsilon_{t2^n+k}(t2^n+k)^j
=
\sigma_t b_t.
\]
Thus the block-endpoint partial sums differ from a fixed initial constant by a signed sum of the block increments \(b_t\).

Let \(q_1,q_2,q_3,\dots\) be a sequence of rational numbers in which every rational occurs infinitely often. We shall choose the block signs in stages. At stage \(\ell\), using only unused block increments, apply Lemma \ref{lem:steer} to choose finitely many further block signs so that the next block-endpoint partial sum lies within \(1/\ell\) of \(q_\ell\). This is possible because every unused tail of the sequence \((|b_t|)\) still tends to \(0\) and has divergent sum.

This constructs a single infinite sequence of signs \(\varepsilon_k\in\{\pm1\}\) and a sequence of block endpoints at which the corresponding partial sums approach
\(
q_\ell
\)
within \(1/\ell\). To see that these block-endpoint partial sums are dense in \(\mathbb R\), let \(r\in\mathbb R\) and let \(\epsilon>0\). Choose a rational \(q\) with \(|q-r|<\epsilon/2\). Since \(q\) occurs infinitely often among the \(q_\ell\), choose \(\ell\) sufficiently large that \(q_\ell=q\) and \(1/\ell<\epsilon/2\). At the corresponding block endpoint,
\[
|S_N-r|\le |S_N-q|+|q-r|<\epsilon.
\]
Hence the partial sums are dense in \(\mathbb R\). Equivalently, for every \(r\in\mathbb R\), some subsequence of the partial sums converges to \(r\).
\end{proof}

\section{Further examples}\label{sec:general}

The proof of Theorem \ref{thm:main-theorem} uses only the fact that the
Thue--Morse block values tend to \(0\), but not so quickly that their absolute
values are summable. This allows the same argument to handle some natural
perturbations of power functions.

\begin{theorem}[Logarithmic perturbations]\label{thm:log-perturbations}
Let \(j>0\) be non-integer and let \(\beta\in\mathbb R\). Then there is a
sequence of signs \(\varepsilon_k\in\{\pm1\}\) such that the partial sums
\[
\sum_{k=2}^{N}\varepsilon_k k^j(\log k)^\beta
\]
are dense in \(\mathbb R\).

Also, if \(m\ge1\) is an integer and \(\beta>0\), then there is a sequence of
signs \(\varepsilon_k\in\{\pm1\}\) such that the partial sums
\[
\sum_{k=2}^{N}\varepsilon_k k^m(\log k)^\beta
\]
are dense in \(\mathbb R\).
\end{theorem}

\begin{proof}
We give the proof for the first statement; the second is similar. Let
\[
f(x)=x^j(\log x)^\beta
\]
and choose \(n=\lceil j\rceil\). Repeated differentiation gives
\[
f^{(n)}(x)=j(j-1)\cdots(j-n+1)x^{j-n}(\log x)^\beta
+O\left(x^{j-n}(\log x)^{\beta-1}\right).
\]
Since \(j\) is non-integer, the leading coefficient is nonzero. Taylor's
theorem and the moment cancellation of the length-\(2^n\) Thue--Morse block
therefore give
\[
(B_n f)(x)=C_{j,n,\beta}x^{j-n}(\log x)^\beta
+O\left(x^{j-n}(\log x)^{\beta-1}\right)
+O\left(x^{j-n-1}(\log x)^\beta\right),
\]
where \(C_{j,n,\beta}\ne0\). Hence the block values tend to \(0\). Since
\(j-n\in(-1,0)\), their absolute values are not summable along the block
endpoints. The same block-sign construction used in the proof of Theorem
\ref{thm:main-theorem} now gives a single choice of signs for which the
corresponding partial sums are dense in \(\mathbb R\).

For the second statement, take
\[
f(x)=x^m(\log x)^\beta
\]
and use the length-\(2^{m+1}\) Thue--Morse block. In this case
\[
f^{(m+1)}(x)=\beta m! x^{-1}(\log x)^{\beta-1}
+O\left(x^{-1}(\log x)^{\beta-2}\right).
\]
For \(\beta>0\), these block values tend to \(0\), while their absolute values
are not summable, since
\[
\sum_{k\ge2}\frac{1}{k(\log k)^{1-\beta}}
\]
diverges. The same steering argument applies.
\end{proof}

\begin{remark}
The lower limit \(k=2\) is used only to avoid the value \(\log 1=0\) when
negative powers of \(\log k\) occur. Changing or deleting finitely many initial
terms has no effect on the density conclusion.
\end{remark}

\section{Further directions}
One may ask for a characterisation of sequences \((a_k)\) for which there is a single choice of signs \(\varepsilon_k\in\{\pm1\}\) such that the partial sums
\(
\sum_{k=1}^{N}\varepsilon_k a_k
\)
are dense in \(\mathbb R\). The present paper gives one mechanism: large terms are first grouped into small non-summable block increments. Other mechanisms may also be possible.

Oscillatory perturbations appear to be a particularly interesting test case. For example, do the signed partial sums of
\(
a_k=k^j\sin k
\)
admit a dense signing? This is not covered directly by our Thue--Morse block argument. It therefore seems to require a different blend of cancellation and equidistribution. 

Perhaps a more natural oscillatory test case is the signed sum
\(
\sum_{k=1}^{N}\sigma_k\sin k.
\)
Here the terms are bounded but do not tend to zero. Since the values \(\sin k\) are dense in \([-1,1]\), there are arbitrarily small available terms, but unlike the block sums in the present paper these small terms are not produced by a fixed cancellation pattern. It would be interesting to know whether there is a
single choice of signs \(\sigma_k\in\{\pm1\}\) for which the partial sums are dense in \(\mathbb R\).

\end{document}